\documentclass{amsart}
\usepackage{a4wide}
\usepackage{amsfonts}
\usepackage{url}
\usepackage{listings}
\usepackage[osf,sc]{mathpazo}
\theoremstyle{plain}
\newtheorem{theorem}{Theorem}
\usepackage{tcolorbox}

\usepackage{breqn}

\theoremstyle{definition}

\allowdisplaybreaks

\subjclass[2020]{11Y65}

\keywords{Rogers-Ramanujan continued fraction, Jacobi's theta function} 

\author{Sumit Kumar Jha}
\address{G-2, C-185, Ramprastha, Ghaziabad, U. P. - 201011}
\email{kumarjha.sumit@research.iiit.ac.in, sumitkumarjha.iiit@gmail.com} 

\begin{document}
	
\title{Two complementary relations for the Rogers-Ramanujan continued fraction}

\begin{abstract}
Let $R(q)$ be the Rogers-Ramanujan continued fraction. We give different proofs of two complementary relations for $R(q)$ given by Ramanujan and proved by Watson and Ramanathan. Our proofs only use product expansions for classical Jacobi theta functions. 
\end{abstract}

\maketitle

\section{Introduction}
The Rogers-Ramanujan continued fraction, denoted by $R(q)$, is defined by
$$
R(q):={\cfrac {q^{1/5}}{1+{\cfrac {q}{1+{\cfrac {q^{2}}{1+{\cfrac {q^{3}}{1+\ddots }}}}}}}}\qquad |q|<1.
$$
We assume $|q|<1$ hereon. Ramanujan \cite{Ram} proved that
\begin{equation}
\label{first}
    R(q) = q^{1/5}\, \frac{\sum_{\lambda=-\infty}^{\infty}(-1)^{\lambda}\, q^{\frac{\lambda}{2}(5\, \lambda+3)}}{\sum_{\lambda=-\infty}^{\infty}(-1)^{\lambda}\, q^{\frac{\lambda}{2}(5\, \lambda+1)}}=q^{1/5}\, \prod _{n=1}^{\infty }{\frac {(1-q^{5n-1})(1-q^{5n-4})}{(1-q^{5n-2})(1-q^{5n-3})}}
\end{equation}
where the second equality is due to the Jacobi's triple product identity \cite{Andrews}. \par 
Ramanujan in his notebooks states the following results:
\begin{theorem}[{\cite[p. 58]{Bruce2}}]\label{thm1}
If $a$ and $b$ are positive and $ab = 1$, then
\begin{equation}
\label{mainone}
\left
\{
\frac{\sqrt{5}+1}{2} + R\left(e^{-2\pi a}\right)
\right
\}
\left
\{
\frac{\sqrt{5}+1}{2} + R\left(e^{-2 \pi b}\right)
\right
\}
=
\frac{5+\sqrt{5}}{2}.
\end{equation}
\end{theorem}
\begin{theorem}[{\cite[p. 91]{Bruce2}}]
\label{thm2}
If $a$ and $b$ are positive and $ab= 1/5$, then
\begin{equation}
\label{maintwo}
\left
\{
\left(\frac{\sqrt{5}+1}{2}\right)^{5} + R^{5}\left(e^{-2\pi a}\right)
\right
\}
\left
\{
\left(\frac{\sqrt{5}+1}{2}\right)^{5} + R^{5}\left(e^{-2 \pi b}\right)
\right
\}
= 5\, \sqrt{5}\,
\left(\frac{\sqrt{5}+1}{2}\right)^{5}.
\end{equation}
\end{theorem}
The above results were proved by Watson \cite{Watson} and Ramanathan \cite{Ramanathan}. We show that these results directly follow from four identities due to Ramanujan in his lost notebook and the product expansion of one of Jacobi's theta functions.
\section{Jacobi's theta functions}
The second Jacobi theta function, denoted by $\vartheta _{2}(z;q)$, and the fourth Jacobi theta function, denoted by $\vartheta _{4}(z;q)$ \cite[p. 166]{Rademacher}, are defined by
\begin{equation}
\vartheta _{2}(z;q)=\sum _{{n=-\infty }}^{\infty }q^{{(n+1/2)^{2}}}\exp((2n+1)iz),
\end{equation}
\begin{equation}
    \vartheta _{4}(z;q)=\sum _{{n=-\infty }}^{\infty }(-1)^{n}q^{{n^{2}}}\exp(2niz).
\end{equation}
Let $q=e^{i\,\pi\, \tau}$. The Jacobi's identity \cite[p. 177]{Rademacher}
$$
\vartheta_{2}\left(\frac{-z}{\tau};\frac{-1}{\tau}\right)=\sqrt{-i\, \tau}\,\, e^{\frac{i\, z^2}{\tau \pi}}\,\, \vartheta_{4}(z,\tau)
$$
gives us the expression
\begin{equation}
\vartheta_{2}(z;\tau)=\frac{1}{\sqrt{-i\, \tau}}\, \sum_{n=-\infty}^{\infty}(-1)^{n}\, \exp\left\{{\frac{\pi}{i\, \tau}\,\left(n-\frac{z}{\pi}\right)^{2}}\right\}.
\end{equation}
In \eqref{first}, we replace $q$ by $q^{2}$ and complete the square in exponents of $q$ to get
\begin{equation}
R(q^{2}) = \frac{\sum_{\lambda=-\infty}^{\infty}(-1)^{\lambda}\, q^{5\, (\lambda+3/10)^{2}}}{\sum_{\lambda=-\infty}^{\infty}(-1)^{\lambda}\, q^{5\, (\lambda+1/10)^{2}}}
=\frac{\vartheta_{2}\left(\frac{3\, \pi}{10};\frac{-1}{5 \,\tau}\right)}{\vartheta_{2}\left(\frac{ \pi}{10};\frac{-1}{5\, \tau}\right)}.
\end{equation}
The above equation was obtained in \cite{Rogers}. By replacing $\tau$ by $-1/\tau$ in the above equation gives us
\begin{equation}
\label{RamTheta}
R\left(e^{-\frac{2\, i\, \pi}{\tau}}\right) = \frac{\vartheta_{2}\left(\frac{3\, \pi}{10};q^{\frac{1}{5}}\right)}{\vartheta_{2}\left(\frac{ \pi}{10};q^{\frac{1}{5}}\right)}.
\end{equation}
We also recall the product expansion for $\theta_{2}$ \cite[p. 171]{Rademacher}
\begin{equation}
    \label{expansion}
    \vartheta_{2}\left(z;q\right)=2\,q^{1/4}\, \cos z\prod\limits_{n=1}^{\infty}{\left(1-q%
^{2n}\right)}{\left(1+2q^{2n}\cos\left(2z\right)+q^{4n}\right)}.
\end{equation}
\section{Product expansions due to Ramanujan}
In the following, let
$$
\alpha = \frac{1-\sqrt{5}}{2}\qquad \text{and}\qquad \beta=\frac{1+\sqrt{5}}{2}.
$$
Ramanujan gave the following product expansions in his lost notebook:
\begin{theorem}[{\cite[p. 21]{Bruce2}}]
Let $t=R(q)$. Then
\begin{align}
\label{one}
 \frac{1}{\sqrt{t}}-\alpha\, \sqrt{t} &= \frac{1}{q^{1/10}}\, \sqrt{\frac{f(-q)}{f(-q^{5})}}\, \prod_{n=1}^{\infty}\frac{1}{1+\alpha\, q^{n/5}+q^{2n/5}},\\
 \label{two}
 \frac{1}{\sqrt{t}}-\beta\, \sqrt{t} &= \frac{1}{q^{1/10}}\, \sqrt{\frac{f(-q)}{f(-q^{5})}}\, \prod_{n=1}^{\infty}\frac{1}{1+\beta\, q^{n/5}+q^{2n/5}},\\
 \label{three}
 \left(\frac{1}{\sqrt{t}}\right)^{5}-\left(\alpha\, \sqrt{t}\right)^{5} &= \frac{1}{q^{1/2}}\, \sqrt{\frac{f(-q)}{f(-q^{5})}}\, \prod_{n=1}^{\infty}\frac{1}{(1+\alpha\, q^{n}+q^{2n})^{5}},\\
 \label{four}
 \left(\frac{1}{\sqrt{t}}\right)^{5}-\left(\beta\, \sqrt{t}\right)^{5} &= \frac{1}{q^{1/2}}\, \sqrt{\frac{f(-q)}{f(-q^{5})}}\, \prod_{n=1}^{\infty}\frac{1}{(1+\beta\, q^{n}+q^{2n})^{5}}.
 \end{align}
Here $f(-q)= \sum_{n=-\infty}^{\infty}(-1)^{n}\, q^{n(3n-1)/2}$ \cite[p. 11]{Bruce2}.
\end{theorem}
\section{Proof of Theorem \ref{thm1}}
\begin{proof}
Replacing $q$ by $q^2$ in \eqref{one} and \eqref{two} and then dividing one by another gives
\begin{equation}
\frac{1-\beta\, R(q^{2})}{1-\alpha R(q^{2})} = \prod_{n=1}^{\infty} \frac{1+\alpha\, q^{2n/5}+q^{4n/5}}{1+\beta\, q^{2n/5}+q^{4n/5}}.
\end{equation}
Now using \eqref{expansion} we have
$$
\frac{\vartheta_{2}\left(\frac{3\, \pi}{10};q^{\frac{1}{5}}\right)}{\vartheta_{2}\left(\frac{ \pi}{10};q^{\frac{1}{5}}\right)}= \frac{\cos(3\, \pi/10)}{\cos(\pi/10)}\,\prod_{n=1}^{\infty} \frac{1+\alpha\, q^{2n/5}+q^{4n/5}}{1+\beta\, q^{2n/5}+q^{4n/5}}.
$$
Thus \eqref{RamTheta} becomes
$$
R\left(e^{-\frac{2\, i\, \pi}{\tau}}\right) = (-\alpha)\, \frac{1-\beta\, R(e^{2\, i \, \pi\, \tau})}{1-\alpha R(e^{2\, i \, \pi\, \tau})}.
$$
Letting $b = -i \tau$ and $a= 1/b = i/\tau$ in the above gives us \eqref{mainone}.
\end{proof}
\section{Proof of Theorem \ref{thm2}}
\begin{proof}
Replacing $q$ by $q^2$ in \eqref{three} and \eqref{four} and then dividing one by another gives
\begin{equation}
\frac{1-\beta^{5}\, R^{5}(q^{2})}{1-\alpha^{5} R^{5}(q^{2})} = \prod_{n=1}^{\infty} \left(\frac{1+\alpha\, q^{2n}+q^{4n}}{1+\beta\, q^{2n}+q^{4n}}\right)^{5}.
\end{equation}
Now using \eqref{expansion} we have
$$
\left(\frac{\vartheta_{2}\left(\frac{3\, \pi}{10};q\right)}{\vartheta_{2}\left(\frac{ \pi}{10};q\right)}\right)^{5}= \frac{\cos^{5}(3\, \pi/10)}{\cos^{5}(\pi/10)}\,\prod_{n=1}^{\infty}
\left(\frac{1+\alpha\, q^{2n}+q^{4n}}{1+\beta\, q^{2n}+q^{4n}}\right)^{5}.
$$
Thus
$$
\left(\frac{\vartheta_{2}\left(\frac{3\, \pi}{10};q\right)}{\vartheta_{2}\left(\frac{ \pi}{10};q\right)}\right)^{5} = (-\alpha^{5})\, \frac{1-\beta^{5}\, R^{5}(q^{2})}{1-\alpha^{5} R^{5}(q^{2})}.
$$
Then replacing $q$ by $q^{1/5}$ in above and using \eqref{RamTheta} gives us
$$
R^{5}\left(e^{-\frac{2\, i\, \pi}{\tau}}\right) = (-\alpha)^{5}\, \frac{1-\beta^{5}\, R^{5}(q^{2/5})}{1-\alpha^{5} R^{5}(q^{2/5})} = (-\alpha)^{5}\, \frac{1-\beta^{5}\, R^{5}(e^{(2\, \pi\, i\, \tau)/5})}{1-\alpha^{5} R^{5}(e^{(2\, \pi\, i\, \tau)/5})}.
$$
Letting $b = \frac{-i \tau}{5}$ and $a= \frac{1}{5b}=i/\tau$ in the above gives us \eqref{maintwo}.
\end{proof}

\end{document}